\input amstex
\documentstyle{amsppt}
\magnification\magstep1         
\baselineskip=10pt
\parindent=0.2truein
\pagewidth{6.6truein} \pageheight{9.0truein} \topmatter \subjclass{Primary:37F45,
Secondary:37F30.}
\endsubjclass

\title  Remarks on the Dynamic of the Ruelle Operator and invariant differentials.
\endtitle
\author Peter M. Makienko
\endauthor
\abstract{Let $ R $ be a rational map. We are interesting in the
dynamic of the Ruelle operator on suitable spaces of
differentials. In particular
 the necessary and sufficient conditions (in terms of convergence of
sequences of measures) of existence of invariant conformal structures on $ J(R) $
are obtained.}
\endabstract
\rightheadtext{Remarks on Ruelle Operator, invariant quadratic differentials and}
\address{Permanent address: Instituto de Matematicas, UNAM at Cuernavaca\newline
Av. de Universidad s/N., Col. Lomas de Chamilpa,\newline C.P. 62210, Cuernavaca,
Morelos, Mexico.}
\endaddress
\thanks{This work has been partially supported by the Russian Fund of Basic
Researches, Grant 99-01-01006.\hfill\newline This work was partially supported by
proyecto {\bf CONACyT} {\it Sistemas Din\'amicos} 27958E from Mexico}
\endthanks

\endtopmatter
\document
\heading{\bf Introduction and main statements}\endheading
\par Let $ R  $ be a rational map. Let $ c \in J(R) $ be a
critical point with infinite forward orbit which does not contain any other critical
point. Then formally we have the following possibilities: \roster \item
$lim_{n\to\infty} \vert (R^n)'(R(c)\vert = 0, $ \item there exists a subsequence $
\{n_i\} $ such that $ lim_{i\to\infty} \vert (R^{n_i})'(R(c)\vert = \infty, $ \item
there exists a subsequence $ \{n_i\} $ such that $ lim_{i\to\infty} \vert
(R^{n_i})'(R(c)\vert = M < \infty $ and $ M \neq 0.$
\endroster

We believe that the first case contains  a contradiction. Because in this situation
the forward orbit of $ c $ must converge to an attractive (superattractive) cycle
and hence $ c \notin J(R).$ This conjecture is true for example for real quadratic
polynomials whose critical point has strictly negative Lyapunov exponent (see
\cite{MS}).

As for the last two cases, the Fatou conjecture claims that  $ R $ is an unstable
map.
With additional conditions placed on the behavior of Poincar{\'e} series (see
definitions below) we know that in this situation $ R $ is an unstable map see
(\cite{Av, Mak2} ) for maps with non-empty Fatou set (\cite{Av}) and for more
general situation (\cite{Mak2}) and (\cite{Lev, Mak1}) for polynomials of degree
two. By using arguments of (\cite{Av, Lev, Mak1, Mak2} we reproduce this result for
any rational map (see theorem 16). How we know this conditions
appear  independently  in works of P. Makienko \cite{Mak1} and G. Levin \cite{Lev}.
 Unfortunately we can not avoid these additional conditions even in the best case
(the Poincar{\'e} series is absolutely convergent) see (\cite{Av, Lev, Mak2} in this
case these conditions are connected with the following conjecture (see below
"Generalized Sullivan Conjecture"):

{\it Except for the Latt{\'e}s maps, there is no measurable invariant integrable
(over $ \overline{{\Bbb C}}$) quadratic differential for rational maps.}

The infinitesimal content of Thurston's Uniqueness Theorem (see (\cite{DH} and also
\cite{MM1}) is the assertion $ R^*(\phi) \neq \phi $ for non-zero integrable
meromorphic quadratic differential $ \phi, $ where $ R^* $ is an "pushforward" or
"Ruelle" operator (see definition below) associated  with the rational map $ R. $ C.
McMullen (\cite{MM2}) proves that $ R $ is the Latt{\'e}s map if and only if there
exists an invariant integrable meromorphic differential. A. Epstein (\cite{E})
extends the results above to: there is no invariant meromorphic differential on $
\overline{{\Bbb C}}$ for any rational map except the Latt{\'e}s map.  P. Makienko
(\cite{Mak2}) shows that with additional assumptions on postcritical set there is no
invariant integrable differential among the augmented meromorphic differentials
except the Latt{\'e}s map. The differential $\phi $ is {\it augmented meromorphic}
if $ \phi = \sum_{i\geq 0}\gamma_i, $ where $ \gamma_i ) $ are meromorphic
integrable differentials with only four simple poles: three of them are $ 0, 1,
\infty $ and series $ \sum\Vert\gamma_i\Vert $is absolutely convergent.

Let us call a measurable differential $ \phi $ as {\it a regular } iff $
\overline{\partial}\phi $ is a finite complex valued measure on $ \Bbb C, $ here
$\overline{\partial}$ in distributions. Note that any augmented differential is a
regular.

In this paper we extend the  result from {\cite{Mak2} and show that with the same
assumption (like in {\cite{Mak2}) on the postcritical set there is no invariant
regular differential except the  Latt{\'e}s map. (Theorem A).


The next two results (see propositions 14-16) give the necessary
and sufficient conditions (in the terms of dynamic of the Ruelle
operator and its modulus) to absence of the invariant conformal
structures on the Julia set and to equality to $ 0 $ the Lebesgue
measure of the Julia set, respectively.

 In our final results (theorems B and C) we reformulate these mentioned
additional conditions in other terms such as the common behavior
of the Ruelle-Poincar{\'e} series (theorem C) and convergence of
the measures's sequence (theorem B).

The main technical idea is:

 Assume that $ J(R) $ supports a non - trivial invariant conformal
structure $ \mu.$ Let $ f_\mu $ be its corresponding quasiconformal map. We find the
conditions which allow us to construct a quasiconformal map $ h $ supported already
on the Fatou set, so that $ h $ and $ f_\mu $ generate the same infinitesimal
deformation of $ R $ (see also {\cite{Mak}}).

 \subheading{Definitions and Main results}
\par  Let $ R $ be a rational map.

Denote by $ P(R) = \overline{\{\cup_c\cup_nR^n(c), \text{ c is a
critical point}\}} $ the postcritical set of $ R. $

\par  Start again with a rational map $ R $ and consider two actions $ R^*_{n,m} $ and
$ {R_*}_{n,m} $ on a function $ \phi $ at point
$ z $ by the formulas
$$\split
 &R^*_{n,m}(\phi) = \sum\phi(J_i)(J_i')^n(\overline{J_i'})^m =
\sum_{y \in R^{-1}(z)}\frac{\phi(y)}{(R'(y))^n(\overline{R'(y)})^m}, \text { and }\\
&{R_*}_{n.m}(\phi) = \phi(R)\cdot (R')^n\cdot\overline{(R')^m},
\endsplit
$$
where $ n $ and $ m $ are integers and $ J_i, i = 1, ..., d $ are branches
of the inverse map $ R^{-1}. $ Then we have
$$ R^*_{n,m}\circ {R_*}_{n,m}(\phi) = deg(R)\cdot\phi.
$$
In other words, the actions above  consist of the natural actions of $ R $ on the
spaces of  forms of type $ \phi(z)Dz^mD\overline{z}^n. $
\proclaim{Definition}
\roster
\item The operator $ R^* = R^*_{2,0} $ is called {\rm the transfer
(pushforward) operator or Ruelle operator} of the rational map $R.$
\item  The
operator $ \vert R^*\vert = R^*_{1, 1} $ is called {\rm the modulus of the Ruelle
operator}
\item The operator  $ B_R = {R_*}_{-1,1} $ is called {\rm the Beltrami
operator} of the rational map $ R. $
\endroster
\endproclaim
 The operators $ R^* $ and $ \vert R^*\vert $ and their  right inverses
 $ R_*(\phi) = \frac{{R_*}_{2,0}(\phi)}{deg(R)} $ and
 $ \vert R_*\vert(\phi) = \frac{\vert {R_*}_{2,0}\vert(\phi)}{deg(R)} =
 \frac{\phi(R)\vert R'\vert^2}{deg(R)} $
map the space $ L_1(\overline{\Bbb C}) $ into itself with the unit norm. Moreover
the operator $ B_R $ maps the space $ L_\infty(\overline{\Bbb C}) $ into itself with
 unit norm.
\proclaim{Generalized Sullivan Conjecture} Except for the Latt{\'e}s map, there are
no
 integrable supported on the Julia set measurable functions  which are invariant with
 respect to any
 of the
 operators enumerated above.
 \endproclaim
\subheading{$L_p$ -topology} Assume $ 1 < p < \infty, $ then we can define the
action of the Ruelle operator $ (R_\ast)_p : L_p(\overline{\Bbb C}) \mapsto
L_p(\overline{\Bbb C})$ by the following way:
$$(R_\ast)_p(\phi) = \frac{1}{d^{\frac{p - 1}{p}}}\sum
\phi(J_i)(J_i^\prime)^{\frac{2}{p}}.
$$
By analogy with the $ L_1 $ case we define the {\it modulus of $ (R_\ast)_p $} by
the following formula:
$$
\vert(R_\ast)_p\vert(\phi) = \frac{1}{d^{\frac{p - 1}{p}}}\sum
\phi(J_i)\vert(J_i^\prime)\vert^{\frac{2}{p}}.
$$

For definition and properties of the Latt{\'e}s maps we are referee to the preprint
of J. Milnor (see \cite{M}). Follow by J. Milnor we say that a Latt{\'e} map $ R $
is {\it flexible} iff $ J(R) $ posses a non trivial fixed point of the Beltrami
operator $ B_R. $

\proclaim{Proposition A} \roster\item A rational map $ R $ is a
flexible Latt{\'e} map if and only if there exists a $ p > 1, $
such that $ (R_*)_p $ has a fixed point in $ L_p(\overline{\Bbb
C}). $ \item A rational map $ R $ is a Latt{\'e} map if and only
if there exists a $ p > 1, $ such that the modulus $
\vert(R_*)_p\vert $ has a fixed point in $ L_p(\overline{\Bbb
C}).$
\endroster
\endproclaim
\demo{Proof} Proof follows from results of A. Zdunik (\cite{Zd})
and mean ergodicity lemma (see below).
\enddemo
 \proclaim{Definition} Let $ R \in  {\Bbb C}P^{2d+1} $ be a rational
map. The component of {\bf J}-stability of $ R $ is the following space.
$$
\aligned qc_J(R) = \bigl\{
&F \in {\Bbb C}P^{2d+1}: \text{ there are neighborhoods } U_R
\text{ and } U_{F} \text{ of }{ \bold J}( R)\,\, \text{ and }
{\bold J}(F),\\
&\text{respectively and a quasiconformal homeomorphism } h_F : U_R \rightarrow
 U_{F} \text{ such that }\\
&F = h_F \circ R \circ {h^{-1}_F}\bigr\}\big/PSL_2(\Bbb C).
\endaligned
$$
\endproclaim
\proclaim{Definition} Let $ R $ be a rational map. Then the space of {\rm invariant
conformal structures} or {\rm invariant line fields} or {\rm Teichmuller space}  $
T(J(R)) $ of $ J(R) $ is the following space
$$
T(J(R)) = \{Fix(B_R)_{\vert L_\infty(J(R))}\}.
$$
where $ Fix(B_R)\subset L_\infty(J(R) $ is the space of invariant
functions for the Beltrami operator $ B_R. $ A result of D.
Sullivan (see \cite{S}) gives a bound of $ 2deg(R) - 2 $ for the
dimension of $ T(J(R)). $

\endproclaim

\proclaim{Definition} Ruelle-Poincar{\'e} series. \roster \item {\rm Backward
Ruelle-Poincar{\'e} series.}
$$ RS(x, R, a) = \sum_{n = 0}^\infty(R^*)^n(\tau_a)(x), $$
where  $ \tau_a(z) = \frac{1}{z - a} $ and $ a \in \overline{\Bbb C} $ is a parameter.
The series
$$ S(x, R) = \sum_{n = 0}^\infty\vert R^*\vert^n({\bold 1}_{\Bbb C})(x)
$$
is called {\rm the Backward Poincar{\'e} series}. \item {\rm Forward
Ruelle-Poincar{\'e} series.}
$$ RP(x,R) = \sum_{n = 0}^\infty \frac{1}{(R^n)'(R(x))}.
$$
The series
$$ P(x,R) = \sum_{n = 0}^\infty \frac{1}{\vert(R^n)'(R(x))\vert}
$$
is called {\rm the forward Poincar{\'e} series.} The series
$$ A(x, R, a) = \sum_{n = 0}\frac{1}{(R^n)'(a)(x - R^n(a))}
$$
is called {\rm the modified Ruelle-Poincar{\'e} series.}
\endroster
\endproclaim
Note that the Ruelle-Poincar{\'e} series generalize  the Poincar{\'e} series
introduced by C. McMullen for rational maps (see \cite{MM}).

\proclaim{Definition} We call a measurable integrable function $ \phi $ {\rm a
regular} iff $ \overline{\partial}\phi $ is a complex-valued finite measure on $
\overline{\Bbb C}, $ here derivative in the sense of distributions.
\endproclaim

\proclaim{Definition} Let $ X $ be the space of rational maps $ R $ that satisfy one
of the following conditions
 \roster \item  the diameters of all
components of $ {\Bbb C} \backslash P(R) $ are uniformly bounded
away from zero or \item $ m(P(R)) = 0, $ where $ m $ is the
Lebesgue measure or \item $ \{P(R)\cap J(R)\} \subset \cup\partial
D,$ where union runs over components $ D \subset F(R). $
\endroster
\par Note that the last condition holds for  maps with completely invariant domain.
\endproclaim
 \proclaim{Theorem A} Let $ R \in X $ be a rational map. Then either $ R $ is a Latt{\'e}s map or
 there is no regular fixed points for $
R^*:L_1(\overline{\Bbb C})\mapsto L_1(\overline{\Bbb C}).$
\endproclaim


\par Now, we give necessary and sufficient conditions for the existence of measurable
invariant conformal structures on the  Julia set in the terms of
especial sequences of measures.

Let $ U $ be a neighborhood of $ J(R). $ We say that $ U $ -is
{\it an essential neighborhood} iff \roster \item $ U $ does not
contain disks centered at all attractive and superattractive
points and \item $ R^{-1}(U) \subset U. $
\endroster

\proclaim{Definition} Let us define the space $ H(U) \subset C({\overline U}), $
where $ C({\overline U}) $ is the space of continuous functions and
\roster
\item $
U $ is an essential neighborhood of $ J(R) $ and
\item $ H(U) $ consists of $ h\in
C({\overline U}) $ such that $ \frac{\partial h}{\partial \overline{z}} $ (in the
sense of distributions) belongs to $ L_\infty(U) $
\item $ H(U) $ inherits the
topology of $ C({\overline U}).$
\endroster
\medskip
\par Measures $ \nu^i_l. $
\par Here we assume that $ 0, 1 $ and $ \infty $ are fixed points for $ R. $
\roster     \item Let $ c_i $ and $ d_i $ be the critical points and critical
values, respectively. Let $ \gamma_{a}(z) = \frac{a(a-1)}{z(z - )(z - a)}. $ Then
define $ \mu_n^i = \frac{\partial}{\partial\overline{z}}((R^*)^n(\gamma_{d_i}(z)) $
(in sense of distributions). We will show below that $  (R^*)^n(\gamma_{d_i}(z))$ $
= $ $\sum_{i=1}^{2deg(R) -2}\sum_{j = 0}^n \alpha^i_j\gamma_{R^j(d_i)}(z) $ and
hence $ \mu_n^i = \sum_{i=1}^{2deg(R) -2}\sum_{j = 0}^n\alpha^i_j\{(R^j(d_i) -
1)\delta_0 - R^j(d_i)\delta_1 + \delta_{R^j(d_i)}\}, $ where $ \delta_a $ denotes
the delta measure with mass at the point $ a. $ \item Define by $ \nu^i_l $ the
average $ \frac{1}{l}\sum_{k = 0}^{l - 1}\mu_k^i. $
\endroster
\endproclaim
\par In general the coefficients $ \alpha_j^i $ in definition above can be expressed
 as a combinations
of  elements of the Cauchy product $ RP(c_i,R,d_i)\otimes
RP(c_i,R) $ of Ruelle-Poincar{\'e} series. \proclaim{Theorem B}
Let $ R \in X $ be a rational map with the simple critical points
and no simple critical relations.
Let $ c_1 \in J(R) $ be a critical point and $
d_1 = R(c_1) $ be its critical value. If there exist an essential neighborhood $ U $
and a sequence of integers $ \{l_k\} $ such that the measures $ \{\nu^1_{l_k}\} $
converge in the $*$-weak topology on $ H(U), $ then the map $ R $ is not
structurally stable ( or is unstable map).
\endproclaim
\proclaim{Corollary B}\roster\item Under assumptions of the
theorem B the space $ T(J(R)) = \emptyset $ if and only if there
exist an essential neighborhood $ U $ and a sequences of integers
$ \{l_k\} $ such that the measures $ \{\nu^i_{l_k}\} $ converge in
the $*$-weak topology on $ H(U) $ for any $ i = 1, ..., 2deg(R) -
2.$ \item If a rational map  $ R \notin X, $  then $ T(J(R)) =
\emptyset $ if and only if there exist an essential neighborhood $
U $ and a sequences of integers $ \{l_k\} $ such that the measures
$ \{\nu^i_{l_k}\} $ converge to  zero in the $*$-weak topology on
$ H(U) $ for any $ i = 1, ..., 2deg(R) - 2.$
\endroster
\endproclaim

Let $ P $ be a structurally stable polynomial. Consider the
following decomposition
$$
\frac{1}{(P^n)^\prime(z)} = \sum \frac{b_i^n}{z - c_i}, $$ where $
c_i $ are critical points of the polynomial $ P^n. $  Let
$$ B_n = \sum_i \vert b^n_i\vert.$$

Now \proclaim{Definition} A structurally stable polynomial $ P $
is {\rm strongly convergent} iff the sums

$$ B_n < C < \infty.
$$
independently on $ n. $
\endproclaim
\proclaim{Remark} We note that, in some sense, any polynomial is
close to be strongly convergent. Indeed the boundedness of $ B_n $
is equivalent to the boundedness of the total variation of the
following  measures
$$ \overline{\partial} \left(\frac{1}{(P^n)^\prime}\right) = \sum b^n_i \delta_{c_i}.
$$
Let  $ h $ be a polynomial then easy calculations show
$$
\sum b^n_i h(c_i) = \int_\gamma \frac{h(z)}{(P^n)^\prime(z)} dz
\to 0 \text{ for } n \to\infty,
$$
here $ \gamma $ is a closed Jordan curve enclosing the closure of
critical points of $ P^n $ for all $ n. $
\endproclaim

\proclaim{Theorem C} Let $ P \in X $ be a strongly convergent
polynomial. Then $ P $ is a hyperbolic polynomial.
\endproclaim

Finally we shortly discuss the following interesting question

How to interpretate the coefficients of a rational map from point
of view of  the dynamic or what is the  dynamical meaning of the
coefficients?

\par The next results clarify the dynamical meaning of the coefficients $ b^n_i $ coming from the decomposition
of $ \frac{1}{(P^n)^\prime(z)} $  and give some formal relations
between Ruelle-Poincar{\'e} series.

\proclaim{Definition}  We denote the Cauchy product of series $ A
$ and $ B $ by $ A\otimes B $. Let us recall that if $ A = \sum_{i
= 1} a_i $ and $ B = \sum_{i = 1} b_i, $ then $ C =  A\otimes B  =
\sum_{i = 1} c_i, $ where $ c_i = \sum_{j = 1}^i a_jb_{i - j}.$
\endproclaim

\proclaim{Proposition C1} If $ P $ is a structurally stable
polynomial, then
$$ B_n = \sum_{c \in Cr(P)}\frac{1}{\vert P^{\prime\prime}(c)\vert}\sum_{j = 1}^{n -1}\frac{1}{\vert (P^{n - j - 1})'(P(c))\vert}\sum_{P^j(y) = c}\frac{1}{\vert(P^j)'(y)\vert^2},
 $$
and hence we have the following formal equality
$$ \sum_{n = 2} B_n = \sum_{c \in Cr(P)}\frac{1}{\vert P^{\prime\prime}(c)\vert}S(c, P)\otimes P(c,
P),
$$

where $ \otimes $ means the Cauchy product of the series.
\endproclaim
\bigskip
\bigskip

 \proclaim{Proposition C2} Let $ R $ be a rational map with simple
critical points. Let $\infty $ be a fixed point for $ R. $ Then
there exist the following formal relations.
$$ \align
RP(a, R) - 1 &= \sum_i\lambda^i
-\sum_i\frac{1}{R''(c_i)}RS(c_i,R,a)\otimes RP(c_i,R),\text{ where
}
\lambda \text{ is the multiplier of } \infty\\
RS(x, R, a) &=A(x, R, a) + \sum_k\frac{1}{R''(c_k)}A(c_k, R,
a)\otimes  RS(x,R,R(c_k)),
\endalign
$$
where $ c_k $ are the critical points of $ R. $
\endproclaim

\proclaim{Corollary C} Let $ P $ be a structurally stable
polynomial. Assume that for any critical point $ c $ there exists
a constant $ M_c $ so that
$$ \frac{1}{\vert (P^n)^\prime(P(c))\vert} \leq \frac{M_c}{n}  \text{ and } \vert R^*\vert^n({\bold 1}_{\Bbb
C})(c) \leq \frac{M_c}{n},
$$
then $ P $ is a strongly convergent polynomial.
\endproclaim

The last section contains a little discussion and open question
with respect to spectrum of the Ruelle operator and related
topics.

This paper is based on the preprint (\cite{Mak1}).

\subheading{Acknowledgement} I would like to thank to IMS at SUNY Stony Brook FIM
ETH at Zurich and IM UNAM at Cuernavaca for its hospitality during the preparation
of this paper.

\heading{\bf Quadratic differentials for rational maps}\endheading
\par Let $ S_R $ be the Riemann surface associated with the action of $ R $
on its Fatou set, then (see \cite{S}) $ S_R $ is a finite union $ \cup_iS_i $ of
punctured torii,  punctured spheres and foliated surfaces.
\par Let $ A(S_R) $ be the space of quadratic holomorphic integrable
differentials on $ S_R. $ If $ S_R = \cup^N_iS_i, $ then $ A(S_R) = A(S_1) \times
...\times A(S_N), $ where $ A(S_i) $ is the space of quadratic holomorphic
integrable differentials on $ S_i. $ \subheading{Quadratic differentials for
foliated surfaces} By results of  (\cite{S}) a foliated surface $ S $ is either
unit disk or a round annulus with marked points, equipped with a group { \bf G}$_D $
of rotations. This group  { \bf G}$_D $ is an everywhere dense subgroup in the group
of all rotations of $ S $ in the topology of uniform convergence on $ S. $ Hence for
any $ z \in S $ the closure of the orbit  { \bf G}$_D(z) $ yields a circle which is
called {\it a leaf of the invariant foliation}. If the leaf $ l $ contains a marked
point $ x, $ then we call $ l $  {\it a critical leaf} and denote it by $ l_x. $
With the exception of a single  case the boundary $ \partial S $ consists of
critical leafs. This exception is the surface corresponding to the full orbit of a
simply connected superattractive periodic component containing only one
 critical point.
In the latter case the surface $ S $ does not contain critical leafs. In this case,
the modulus of $ S $ is not defined (see {\cite{S}} for
 details).
 \par Any quadratic absolutely integrable
 holomorphic differential $ \phi  $ must be invariant under
 the action of the group { \bf G}$_D $ for foliated surfaces.  Hence  $ \phi = 0 $ is the
only absolutely
 integrable holomorphic
 function for $ S $ with undefined modulus, and therefore we have in this case
 $ A(S) = \{0\}.$
\par
 After removing the critical leaves from $ S $
 we obtain the collection $ \cup S_i\cup D $ of  rings $ S_i $ and disk $ D $
 (in the case of Siegel disks). We call this decomposition as {\it critical
 decomposition.} For this decomposition we have
 $$
 \phi_{\vert S_i} =  h_i(z)\cdot dz^2, \text{ and } \phi_D = h_0\cdot dz^2
  $$
  where $ h_0 $ and $ h_i $ are  holomorphic absolutely integrable functions on $ D $ and
  $ S_i $ respectively.
The   calculations show that $ h_i(z) = \frac{c_i}{z^2} $ and $ h_0 = 0, $ where $
c_i $ are arbitrary constants.
  From the discussion above we conclude that for a ring with $ k $ critical
  leaves (two of which represent the boundary of $ S $) the dimension
  $ dim(A(S)) = k - 1. $
  \par Now let $ S $ be a ring with critical decomposition $ \cup_{i = 1}^k S_i$ and
  $ \phi \in A(S) $  a differential, then $ \Vert\phi\Vert = 4\pi\sum_i\vert
  c_i\vert mod(S_i), $ where $ \phi = \sum_i\phi_{\vert S_i} = \sum_i
  {\frac{c_i}{z^2}}_{\vert S_i} $ and $ mod(S_i) $ is the modulus (or the extremal
  length of the family of curves connecting the boundary component of $ S_i$) of the ring
  $ S_i. $
  \par We always assume here that the hyperbolic metric $ \lambda $ on the foliated
  ring $ S $ is the collection of complete hyperbolic metrics $ \lambda_i $ on
the   components of the critical decomposition of $ S. $ For example if $ \cup_i S_i
$ is
  the critical decomposition of $ S, $ then the space $ HD(S ) $ of harmonic
  differentials on $ S $ consists of the elements
  $$ \lambda^{-2}\overline{\phi} = \sum_i\frac{\overline{c_i}\lambda_i^{-2}}{\overline{z^2}}
  ,
   $$
  where $ \phi \in A(S). $

In fact, the space $ HD(S) $ is isomorphic to the dual space $ A(S(R))^* $ by way of
the Petersen inner product
$$
 <\psi, \phi>  = \iint_{S(R)} \lambda^{-2}\overline{\psi}\phi.
$$

  \par The space of Teichmuller differentials $ td(S) $ consists of the elements
  $ \phi =
  \sum_i {c_i\frac{z}{\overline{z}}\frac{d\overline{z}}{dz}}_{\vert S_i},
  $ where $ \cup_i S_i $ is the critical decomposition of $ S. $
\par
Denote by $ \Omega (R) $ the set
$$ \overline{\Bbb C} \backslash \overline{\{\cup_nR^{-n}(P(R))\}},
$$
 then $ R $ acts on $
\Omega (R) $ as an unbranched autocovering.

\par Now let $ Y \subset \overline{\Bbb C} $ be an open subset. Then, as
above, $ A(Y) $ denotes the space of holomorphic functions on $ Y $
absolutely integrable
over $ Y $ and $ B(Y) $ consists of holomorphic functions
$ \phi $ on $ Y $ with the following norm
$$ \Vert\phi\Vert = sup_{z \in Y}\vert\Lambda_Y^{-2}\phi\vert, $$
where $ \Lambda_Y $ is a  metric so that the its restriction to any component $ D
\subset Y $ satisfies
$$ {\Lambda_Y}_{\vert D} = \lambda_D,
$$
where $ \lambda_D $ is Poincar{\'e} metric on $ D.$ \proclaim{Lemma(Bers Duality
Theorem)}
Let $ Y \subset \overline{\Bbb C} $ be an open subset. Then the spaces $ A(Y) $ and
$ B(Y) $ are dual by the Peterson scalar product
$$ \iint_Y(\Lambda_Y)^{-2}\phi\overline{\psi}.
$$
\endproclaim
\demo{Proof} See \cite{Kra}. \enddemo \subheading{Poincar{\'e} $\Theta-$operator for
rational maps} We construct this operator in a way which is reminiscent of the
construction  in the case of a  Kleinian group.
\par 1). Let $ D \in \Omega (R) $  correspond to an attractive periodic domain.
Let $ S_D \in S(R) $ be corresponding riemann surface. Then the projection $ P:L(D)
\mapsto S_D $ is a holomorphic unbranched covering. Let  $ P^*: A(L(D))\mapsto
A(S_D) $ be the push-forward operator: locally in charts $ P^*(\phi) = \sum (\phi
dz^2)\circ r_i, $ where summation taken over all branches $ r_i $ of $ P $. Then we
call $ P^* -$ {\it Poincar{\'e} operator } for the attractive domain $ D $ and
denote it by $ \Theta_{L(D)}. $

2) The case of parabolic domains $ D $ is similar to that of an attractive
domain.\newline

3) Foliated case. This case corresponds to non-discrete groups. We need additional
information pertaining to the foliated case. Let us start with a simple lemma about
the Ruelle operator. \proclaim{Lemma 1} Let $ R $ be a rational map,  $ Y \subset
\overline{\Bbb C} $ a positive Lebesgue measure subset which is completely invariant
under the action of $ R. $ Then the following is true. \roster

\item $ R^* : L_1(Y) \rightarrow L_1(Y) $ is a linear surjection with unit norm. The
operator
$$ R_*(\phi)= \frac{\phi(R)(R')^2}{deg(R)}
$$ is an isometric inclusion "into" and $ R^*\circ R_* = I, $ where $ I $ is the identity
operator.
\item The Beltrami operator
$$
 B_R(\phi) = \phi(R)\frac{\overline{R'}}{R'}: L_\infty(Y) \rightarrow L_\infty(Y)
 $$
is the dual operator to $ R^*. $  The operator $ B_R $ is an isometric inclusion.
\item If $ Y $ is an open set, then $ R^*:A(Y) \rightarrow A(Y) $ is a surjection of
unit norm and $ R_* $ maps $ A(Y) $ into itself as well.
\endroster
\endproclaim
\demo{Proof} All items are immediate consequences of the definition  of the operators.
\enddemo
\proclaim{Remark 2} Suppose $ R: X \rightarrow Y $ is a branched covering with  a
rational map $ R $ and two domains $ X,Y \subset \overline{\Bbb C}. $ Then $ R^*:
A(X) \rightarrow A(Y) $ is the Poincar{\'e} operator of the covering $ R. $
\endproclaim

\par Now let us move onto the foliated case. For simplicity let $ D $ be
an invariant either superattractive domain or Siegel disk or Herman ring. Our aim is
to prove the following theorem. \proclaim{Theorem 3} Let $ D \subset F(R) $ be an
invariant domain corresponding to the foliated case and $ X = D \backslash P(R). $
Let $ S $ be the foliated surface associated with $ D .$ Then there exists a
continuous linear projection $ P^*: A(X)\mapsto A(S). $ The dual operator $ P_*:
HD(S)\mapsto HD(X) $ is an injection.
\endproclaim

 \demo{Proof}
 Let $ Y \subset X $ be a component of critical decomposition and $ O(Y) $ be its
full orbit, then  the projection $ P:O(Y)\mapsto Y $ allows to construct the
push-forward and pull-back maps. But in this situation the push-forward map maps $
A(O(Y)) $ onto the space $A(Y) $ of all integrable holomorphic quadratic
differentials on $ Y.$ We need to consider average with respect to the group $
G(D).$
\enddemo
We claim

\proclaim{Claim} There exists a continuous linear projection $ \pi:A(X)\mapsto A(S).
$ The dual map $ \pi_* $ is an injection from $ HD(S) $ into the space of Beltrami
differentials on $ X. $
\endproclaim

Firstly assume that claim is proved, then we finish our proposition by setting $ P^*
= \pi \circ \{\text{push-forward map}\} $ and $ P_* = \{\text{pull-back map}\}\circ
\pi_*.$

In what follows we need  the following basic facts about non-expansive operators and
rational quadratic differentials.

\proclaim{Mean ergodicity lemma} Let $ T $ be a non-expansive ($\Vert T\Vert \leq
1$) linear endomorphism of a Banach space $ B $ and let  $ \phi \in B $ be any
element. \roster \item Assume that for the {\rm Cesaro average} $A_N(T,\phi) =
\frac{1}{N}\sum_{i =0}^{N - 1}T^i(\phi) $ there  exists a subsequence $ \{n_i\} $
such that $ A_{n_i}(T,\phi) $ weakly converges to an element $ f \in B. $ Then $ f $
is a fixed point for $ T $ and $ A_N(\phi) $ converges to $ f $ strongly (i.e. in
norm). If $ f = 0 $ then $ \phi \in \overline{(I - T)(B)} $ and vice versa i.e. if $
\phi \in \overline{(I - T)(B)}, $ then $ A_n(T,\phi) $ tends to zero with respect to
the norm. \item A linear continuous operator $ T $ on a norm space $ B $ is called
{\rm mean ergodic} if and only if the Cesaro average $ A_N(T,\phi)$ converges with
respect to the  norm for any element $ \phi \in B.$ In this case $ B = Fix(T)\times
\overline{(I - T)(B)} $ and $ A_n(T) $ converges (in the strong topology) to the
continuous projection $ \pi: B \rightarrow Fix(T), $ where $ Fix(T) $ is the space
of fixed elements for $ T.$ \item A linear non-expansive endomorphism  $ T $ of a
Banach space $ B $ is mean ergodic if and only if for any $ \omega \in Fix(T^*) $
there exists $ \phi \in Fix(T) $ so that $ \omega(phi) \neq 0, $ where $ Fix (T^*) $
is the space of fixed elements for the dual operator $ T^*.$
\endroster
\endproclaim
\demo{Proof}  See  text of U. Krengel (\cite{Kren}), theorem 1.1, page 72 and the
theorem 1.3, page 73.
\enddemo
\proclaim{Lemma(Bers Density Theorem)} Let $ C $ be a closed subset of $
\overline{\Bbb C} $ and $ C_0 $ a dense subset of $ C. $ If $ Y $ is the complement
to $ C $ in $ \overline{\Bbb C}, $ then let $ A(Y) \subset L_1(\overline{\Bbb C}) $
be subspace of the functions holomorphic on $ Y. $ Let $ R(C_0) \subset
L_1(\overline{\Bbb C}) $ be subspace of the rational functions holomorphic outside
of $ C_0.$ Then $ R(C_0) $ is everywhere dense subspace in $ A(Y) $ in the
$L_1-$norm.
\endproclaim

\demo{Proof} See the text  of F.P. Gardiner and N. Lakic {\cite{GL}).
\enddemo
Let $ Z $ be the surface $ Y $ equipped with the group $ G(D).$ Any element $ h \in
A(Z) $ is invariant with respect to the group of all rotation of $ Y $ and vise
versa. Then the average with respect to any dense subgroup $ G \in S^1 $ must
produces elements of $ A(Z).$

Let us consider $ g(z) = \rho z, \rho = \exp(2\pi i \alpha), $ with an irrational $
\alpha.$ Then $ g $ defines an non-expansive operator $ T : A(Y)\mapsto A(Y), $ by
the formula $ T(h)(z) = h(\rho z)\rho^2 $ and the space $ Fix(T) = A(Z).$

Our aim is to show that $ T : A(Y)\mapsto A(Y) $ is mean ergodic.

Let $ W \subset A(Y) $ be the linear span of quadratic integrable rational
differentials on $ \overline{\Bbb C} $ with poles in $ \overline{\Bbb C}\backslash
\overline{Y}. $ Then by the Bers's density theorem $ W $ is everywhere dense subset
of $ A(Y).$ If the Cesaro averages $ A_n(T,h) $ converges in the strong topology for
any element $ h \in W, $ then $ T $ is mean ergodic such as the Cesaro averages form
equicontinuous family of operators.

Hence by above it is enough to show the strong convergence of $ A_n(T,h) $ for any $
h\in W. $ Let $ h \in Y $ be an element, then by assumptions $ h $ is holomorphic in
a ring $ E $ which compactly contains the ring $ Y. $ Let $ h(z) =
\sum_{-\infty}^{+\infty} a_iz^i$ be the  Laurent series of $ h, $ then the
decomposition of $ A_n(T,h) $ is as follows
$$
\sum_{-\infty}^{+\infty} a_iz^i\frac{\sum_{l = 0}^{n - 1}\rho^{l(i+2)}}{n}
$$
Such as $ \rho = \exp(2\pi i\alpha) \neq 1 $ and the series $
\sum_{-\infty}^{+\infty} a_iz^i $ converges absolutely on $ \overline{Y}, $ then the
series for $ A_n(T,h) $ converges absolutely too on $ \overline{Y}. $ Hence the
calculations by elements show that $ \lim_{n\to\infty}\frac{\sum_{l = 0}^{n -
1}\lambda^{l(i+2)}}{n} = 0, $ for any $ i \neq -2. $ as result we have
$$
\lim_{n\to\infty} A_n(T,h)  = \frac{a_{-2}}{z^2}
$$
The uniform convergence in limit above  on $ \overline{Y} $ implies the strong
convergence in $ A(Y).$

The dual map $ \pi_* : HD(Z) \mapsto HD(Y) $ is simply inclusion so that the
following is true
$$ \iint_Y \pi(h)\mu = \iint_Y h \pi_*(\mu),
$$
for any $ h \in A(Y) $ and $ \mu \in HD(Y).$ Hence the claim and the theorem are
proved.

\par Finally we define
$ \Theta(R): A(\Omega)\rightarrow A(S_R) $ by
$$ \Theta(R)(\phi) = \left(\Theta_{O(D_1)},..., \Theta_{O(D_k)}\right),
$$
where $ D_i \subset F(R) $ are periodic components. If $ D_i $ is a superattractive
or Siegel or Herman, then $\Theta_{O(D_i)} = P^*.$ Here $ O(D_i) $ means the full
orbit of $ D_i. $

\heading{\bf The Space $A(R)$}\endheading
\par Now again consider the space $ A(\Omega). $ Note that any function
of the form
$$ \gamma_a(z) = \frac{a(a - 1)}{z(z - 1)(z - a)} \text{ for } a \in
\overline{\Bbb C} \backslash \Omega $$
belongs to $ A(\Omega).$  Let us introduce the subspace $ A(R) \subset
A(\Omega) $ as follows.
Let $ S $ be the set
$$ \{\cup_{i}\{O(c_i)\}\cup \{\{O(0, 1, \infty)\}
\backslash \{0, 1, \infty\}\}, $$ where $ c_i $ are the critical points . Then we
set
$$ A(R) = \text{linear span}\{\gamma_a(z), a \in S\} $$
This space $ A(R) $ is a linear space and we introduce on $ A(R) $ two different
topologies though  the  norms $ \vert\cdot\vert_1 = \int_\Omega \vert\cdot\vert $
and $ \vert\cdot\vert_2 = \int_{{\bold J}(R)} \vert\cdot\vert $. Denote by $ A_i $
the spaces $ \{A(R), \vert\cdot\vert_i\} $, respectively. \proclaim{Remark 4} The
space $ A(R) $ serves as a kind of connection between the spaces $ L_1(\Omega) $ and
$ L_1(J(R)). $ The comparison of the $\Vert\cdot\Vert_1 $ and $ \Vert\cdot\Vert_2 $
topologies is the basis for our discussion below.
\endproclaim
\proclaim{Lemma 5} The operators $ R^* $ and $  R_* $ are continuous endomorphisms
of $ A_1 $ and $ A_2. $
\endproclaim
\demo{Proof} It is sufficient to show that for any $ \phi \in A(R) $ the functions $
R^*(\phi) $ and $R_*(\phi) $ belong to $ A(R). $
\par Let $ \phi = \gamma_a. $ Then $ R^*(\phi) $ and $R_*(\phi) $ are holomorphic
everywhere except a finite number of points belonging to the set $ S $ and hence are
rational functions holomorphic on $\Omega. $ Moreover both $ R^*(\phi) $ and
$R_*(\phi) $ are integrable over $ \overline{\Bbb C} $ and hence belong to $ A(R). $
This proves the lemma.
\enddemo

\proclaim{Lemma 6} Let $ L $ be a continuous functional on $ A_1 $ invariant under
the action of $ R^\ast $ (i.e. $L((R^\ast)(\phi)) = L(\phi)$).Then $ L(\phi) =
\iint\limits_\Omega \lambda^{-2}{\overline{\psi}}\phi, $ where $ \lambda $ is
hyperbolic metric on $ \Omega $ and $ \psi \in B(\Omega).$
\endproclaim
\demo{Proof} The Bers density and Bers duality theorems prove  this lemma.
\enddemo

Now we  need the following basic facts from the Ergodic Theory:

Let $ T: X \mapsto X $ be measurable non-singular map of a measurable subset $ X
\subset {\Bbb C}, $ with respect to the Lebesgue measure $ m. $

Let $ C\cup D = X $ be the Hopf decomposition of $ X $ onto conservative part $ C $
and dissipative part $ D. $ Let us recall that a measurable subset $ W\subset X $ is
{\it wandering} iff $ m(T^{-n}(W)\cap T^{-k}(W)) = 0, $ for any $ k \neq n. $ Then $
D = \{\cup W, W \subset X \text{ is wandering }\}, $ and $ C = X \backslash D. $
Since $ T^{-1}(W) $ is wandering for wandering $ W, $ we have that $ T^{-1}(D)
\subset D $ modulo the Lebesgue measure. Hence $ T(C) \subset C $ modulo the
Lebesgue measure.

 The map $ T: X \mapsto X $ is called {\it
conservative} if $ m(D) = 0.$ If $ m(C) > 0, $ then from above we have that $
T_{\vert C}:C \mapsto C $ is a conservative map.

\proclaim{Poincar{\'e} Recurrence Theorem} Suppose $ T: X\mapsto X $ is a
conservative, non-singular map. If $ (Z, d) $ is a separable metric space, and $ f :
X\mapsto Z $ is a measurable map, then
$$
\lim\inf_{n\to\infty} d(f(x), f(T^n(x))) = 0 \text{ for almost every } x \in X.
$$
\endproclaim
\demo{Proof} See text of J. Aaronson {\cite{\bf Aa}} page 17.\enddemo

\proclaim{Definition}A rational map $ R $ is {\rm ergodic} if any measurable set $ A
$ satisfying $ R^{-1}(A) = A $ has zero or full measure in the sphere. \endproclaim

Then  the Poincar{\'e} Recurrence Theorem above and a result of
Lyubich {\cite L}  ( or see theorem "Attracting or ergodic" p.42
in the book  of  C. McMullen (\cite{\bf MM1}) ) give the following
alternative.
 \proclaim{Lemma 7} Let $ R $ be any rational map.
Let $ C $ be the conservative part of the Hopf decomposition of $J(R).$ Assume that
the Lebesgue measure $ m(C) > 0, $ then \roster\item $ C \subset P(R) $ modulo the
Lebesgue measure, or \item the Julia set is equal to the whole Riemann sphere and
the action of $ R $ is ergodic.
\endroster
\endproclaim

\demo{Proof} Let $ B = C \backslash P(R). $ Assume $ m(B) > 0. $ Then by the
Poincar{\'e} Recurrence Theorem with $ (Z, d) = (\overline{\Bbb C}, d), $ where $ d
$ is the spherical metric and $ f = id, $ we obtain that
$$
\lim\sup d(R^n(x), P(R)) > 0,
$$
for almost every $ x \in B.$ Hence by the arguments of theorem of
 McMullen which are mention above the map $ R $ is
 ergodic. We complete this lemma.
\enddemo

 Now we ready to prove the theorem A.
 \proclaim{Theorem A}Let $ R \in X $ be a rational map. Then either $ R $ is a Latt{\'e}s map or
 there is no regular fixed points for $
R^*:L_1(\overline{\Bbb C})\mapsto L_1(\overline{\Bbb C}).$
\endproclaim
\demo{Proof} Assume $ \phi\neq 0 \in L_1(\overline{\Bbb C}) $ is a
fixed regular point for $R^*.$ Then the dual operator $ B_R:
L_\infty(\overline{\Bbb C})\mapsto L_\infty(\overline{\Bbb C}) $
has a fixed point $ \mu $ so that $ \iint_{\overline{\Bbb
C}}\mu\phi \neq 0. $  The lemma 7 above and theorem 3.17 (Toral or
attracting) in \cite{\bf MM1} imply either $ R $ is a flexible
Latt{\'e} map or the complement to the postcritical set $ \left
\{\overline{\Bbb C} \backslash P(R)\right \} \subset D, $ where $
D $ is the dissipative set for the map $ R. $

 Then $ \vert\phi\vert = 0 $ almost everywhere on $ D $ and hence supporter
 of the complex valued measure $ \overline{\partial}\phi $ belong to $ P(R). $

Now we can assume that $ P(R) $ is a compact subset of the plane. Otherwise  if $
P(R) $ is unbounded let $ h(z) $ be  a Mobius map which maps $P(R) $ into the plane,
then the differential $ \phi_1 = \phi(h)(h^\prime)^2 $ is invariant for the map $
R_1 = h\circ R \circ h^{-1} $ and easy calculations show that $ \overline{\partial}
\phi_1 $ is a complex valued finite measure.

Let $ V(z) = \iint_{\Bbb C} \frac{\overline{\partial}\phi(\xi)}{\xi - z} $ be the
Cauchy transform of the measure $ \overline{\partial}\phi. $ Easy calculations show
that in distributions $ \overline{\partial}V = \overline{\partial}\phi $ and $ V(z)
$ is holomorphic out of $ P(R), $ and $ V(z) \mapsto 0, $ as $ z\mapsto\infty. $
Then we claim: \proclaim{Claim} $ \phi = V $ almost everywhere on the plane.
\endproclaim
\demo{Proof of the claim} Indeed $ \overline{\partial}(\phi - V) = 0 $ in
distribution, hence the function $ \phi - V  $ is holomorphic on the plane. Besides
$ \phi  = 0 $ out of $ P(R) $ and $ \phi(z) - V(z) \mapsto 0, $ as $ z\mapsto\infty$
hence $ \phi = V $ almost everywhere.
\enddemo

To finish the proof of theorem A we need

\proclaim{Proposition 8} Let $ Z $ be a compact subset of the
plane with empty interior and $\omega \neq 0 $ be a complex-valued
finite measure on $ Z. $ Let $ V(z) = \iint_{\Bbb
C}\frac{\omega(\xi)}{\xi - z} $ be Cauchy transform. Then the
function $ V(z) $ is
 not
 identically zero on $ Y = {\Bbb C} \backslash Z $ in the following cases
\roster \item the set $ Z $ has zero Lebesgue measure,
 \item if the diameters of the
components of $ {\Bbb C} \backslash Z $ are uniformly bounded from below away from
zero or \item If $ O_j $ denote the components of $ Y, $ then $ Z \in
\cup_j{\partial O_j}.$
\endroster
\endproclaim
\demo{Proof} The first is evident.

2) Assume that $ V = 0 $ identically outside  of $ Z.$ Let $ R(Z) \subset C(Z) $
denote the algebra of all uniform limits of rational functions with poles outside of
$ Z $ in the $\sup-$topology. Here $ C(Z) $ denotes as usual the space of all
continuous functions on $ Z $ with the $ \sup-$norm. Then the measure $ \omega $
induces a linear functional on $ R(Z).$ Items (2) and (3) are based on the
generalized Mergelyan theorem (see \cite{\bf Gam} thm. 10.4) which states that {\it
If diameters of all components of $ {\Bbb C} \backslash Z $ are bounded uniformly
from below away from 0, then every continuous function holomorphic on the interior
of $ Z $ belongs to $ R(Z). $}

Let us show that $ \omega $ annihilates the space $ R(Z). $ Indeed, let $ r(z) \in
R(Z) $ be a rational function and $ \gamma $ a curve enclosing $ Z $ close enough to
$ Z $ such that $ r(z) $ does not have poles in the interior of $ \gamma.$ Then
since  $ l = 0 $ outside of $ Z $ we need only apply Fubini's theorem:
$$
\int r(z)d\omega(z) = \int d\omega(z)\frac{1}{2\pi
i}\int_{\gamma}\frac{r(\xi)d\xi}{\xi - z} = \frac{1}{2\pi
i}\int_{\gamma}r(\xi)d\xi\int\frac{d\omega(z)}{\xi - z} =
\frac{1}{2\pi i}\int_{\gamma}r(\xi)V(\xi)d\xi = 0.
$$
By the generalized Mergelyan theorem, we have $ R(Z) = C(Z) $ and $ \omega = 0, $
contradiction.

Now let us check (3). We {\bf claim} {\it that $ V = 0 $ almost everywhere on} $
\cup_i\partial O_i. $

{\it Proof of the claim.}  Let $ E \subset \cup_i\partial O_i $ be any measurable
subset with positive Lebesgue measure. Then the function $ F_E(z) =
\iint_E\frac{dm(\xi)}{\xi - z} $ is continuous on $ {\Bbb C} \backslash \cup_i O_i $
and is holomorphic on the interior of $ {\Bbb C} \backslash  \cup_i O_i. $ Again, by
the  generalized Mergelyan theorem  $ F_E(z) $ can be approximated on $ {\Bbb C}
\backslash \cup_i O_i $ by functions from $ R({\Bbb C} \backslash  \cup_i O_i); $
and hence by the above arguments and by our hypotheses we have $ \int
F_E(z)d\omega(z) = 0. $ But another application of Fubini's theorem gives
$$ 0 = \int F_E(z)d\omega(z) = \int d\omega(z)\iint_E\frac{dm(\xi)}{\xi - z} =
\iint_E dm(\xi)\int d\omega(z)\frac{1}{\xi - z} = \iint_E
V(\xi)dm(\xi).
$$
Hence for any measurable $ E \subset \cup_i\partial O_i, $ we have
$ \iint_E V(z) = 0. $ This proves the claim.
\par Now for any component $ O \in Y $ and  any
measurable $ E \subset \partial O $ we have $ \iint_E V(z) = 0. $
By assumption $ V = 0 $ almost everywhere on $ \Bbb C,$
contradiction with $ \omega \neq 0. $ This finishes the proof of
the proposition.
\enddemo

To finish the theorem A in we only need to note that under our assumption the
postcritical set $ P(R) $ has empty interior.

\enddemo

 \heading{\bf  Bers Isomorphism} \endheading
\par Here we reproduce the Bers construction for  Beltrami
differentials and Eichler cohomology  with corrections (which are often obvious) for
the rational maps.
\par Consider the Beltrami action of $ R $ on the space $
L_\infty(\overline{\Bbb C}) $ i.e.
$$ B_R(\phi)(z) = \phi(R)(z)\frac{\overline{R'(z)}}{R'(z)}. $$
Thus the subspace $ Fix $ of fixed points for $ B_R $ in $ L_\infty(\overline{\Bbb
C}) $ is in fact the space of  invariant Beltrami differentials for $ R. $ The  unit
ball in this space describes all quasiconformal deformations of $ R.$

\par Now normalize $ R $ so that 0, 1, $\infty $ are fixed points for $R.$ Let $ K(R) $
be the component of the subset of rational maps in $ {\Bbb C}P^{2d + 1} $ fixing the
points 0, 1 and $ \infty $ containing $ R. $

\par Let $ \mu \in Fix(R); $ then for any $ \lambda $ with
$ \vert\lambda\vert < \frac{1}{\Vert\mu\Vert}, $ the element $ \mu_\lambda =
\lambda\mu \in B. $ Let $ f_\lambda $ be the family of qc-maps corresponding to the
Beltrami differentials $ \mu_\lambda $ with $ f_\lambda(0, 1, \infty) = (0, 1,
\infty).$ Then the map
$$ \lambda \rightarrow R_\lambda = f_\lambda\circ R\circ f_\lambda^{-1} \in K(R)
$$
is a conformal map. Let $ R_\lambda(z) = R(z) + \lambda G_\mu(z) + ... $.
Differentiation with respect to $ \lambda $ at the point $ \lambda = 0$ gives the
following equation
$$ F_\mu(R(z)) - R'(z)F_\mu(z) = G_\mu(z),
$$
where $ F_\mu(z) = \frac{\partial}{\partial\lambda}f_\lambda(z)_{\vert\lambda = 0} $
and $ G_\mu(z) = \frac{\partial}{\partial\lambda}{R_\lambda(z)}_{|\lambda = 0} \in
H^1(R). $
\par By the theory of $qc$ maps (see for example \cite{Krush}). For
any $ \mu \in L_\infty(\Bbb C) $ and $ t $ with $\vert t\vert < \epsilon $ and $
\epsilon$ small, there exists the following formula for the qc-map $ f_{t\mu} $
fixing $ 0, 1, \infty $:
$$
f_{t\mu}(z) = z -\frac{z(z - 1)}{\pi}\iint_{\Bbb C}\frac{t\mu}{\xi(\xi - 1)(\xi -
z)} + \vert t\vert O(C(\epsilon, R)\Vert\mu\Vert_\infty^2),
$$
where $ \vert z\vert < R $ and $ C(\epsilon, R) $ is a constant which does not
depend on $ \mu.$ In addition,
$$
F_\mu(z) = \frac{\partial f_\lambda}{\partial\lambda}_{\vert\lambda = 0} = -
\frac{z(z - 1)}{\pi}\iint_{\Bbb C}\frac{\mu}{\xi(\xi - 1)(\xi - z)}.
$$
\par By $ H^1(R) $ we mean the complex tangent space to $ K(R) $ at the point $ R.$
Then $ H^1(R) $ may be described  as follows: if $ R(z) = z\frac{P_0}{Q_0}, $ then
$$ H^1(R) = \{z\frac{PQ_0 - QP_0}{Q^2_0}, \text{ where } Q(1) =
P(1), deg(Q) \leq deg(R), deg(P) \leq \deg(R) - 1\},
$$
In the above $ P, Q $ are polynomials and $ dim(H^1(R)) = 2d - 2.$ \proclaim{Remark
9} We use the notation $ H^1(R) $ for the following reasons:
 \roster \item The Weyl
cohomology' construction for the action of $ R $ (by the formula $ {\tilde R}(f) =
\frac{f(R)}{R'}$) on the space  of all rational functions gives the space $ H $
which is isomorphic to the tangent space to $ {\Bbb C}P^\infty $ at $ R $ (up to
normalization). More precisely $ H $ is equivalent to the direct limit
$$ H^1(R)\overset{j_1}\to\longrightarrow H^1(R^2)\overset{j_2}\to \longrightarrow
H^1(R^3) ..., $$ where $ j_i $ are equivalent to the action $ \tilde R. $
\item This construction for a Kleinian group gives Eichler cohomology.
\endroster
\endproclaim

\par Hence we can define a linear map $ \beta: Fix(R)\rightarrow H^1(R) $ by the formula
$$
 \beta(\mu) = F_\mu(R(z)) - R'(z)F_\mu(z).
$$
In analogy with Kleinian groups we  call the map $ \beta$ the {\it Bers map }(see
for example \cite{Kra}).
\par Let $ A(S(R)) $ be the space of holomorphic integrable quadratic  differentials on the
disconnected surface $ S(R).$ Let $ HD(S(R)) $ be the space of harmonic
differentials on $ S(R)) $: these are differentials which in a local charts have the
form $ \alpha = \frac{\overline{\phi} d z^2}{\rho^2\vert dz\vert^2}, $ where $ \phi
d z^2 \in A(S(R)) $ and $ \rho\vert d z\vert $ is the Poincar{\'e} metric. In fact,
the space $ HD(S(R)) $ is isomorphic to the dual space $ A(S(R))^* $ by way of the
Petersen inner product
$$
 <\psi, \phi>  = \iint_{S(R)} \rho^{-2}\overline{\psi}\phi.
$$

Let $ \Theta^*:A^*(S_R) \rightarrow A^*(\Omega) $ be the dual operator. Then the
image $ HD(R) = \Theta^*(A^*(S_R))$
is called
{\it the space of harmonic differentials} and $ \dim(HD(R)) = \dim(A^*(S_R)) =
\dim(A(S_R)).$

 By duality we
have
$$
\iint_{\Omega(R)}\Theta^*(\alpha)\phi = \iint_{S(R)}\alpha\Theta(\phi),
$$
for any $ \alpha \in HD(S(R)) $ and $ \phi \in A(\Omega(R)).$ Thus the element $
\beta = \Theta^*(\alpha) $ presents the trivial functional on the space $
A(\Omega(R)) $ if and only if $ \alpha = 0.$

\par Let us recall that $ T(J(R)) = Fix(B_R)_{\vert J(R)} $ is the
space of invariant Beltrami differentials
supported on the Julia set.

The following theorem is proved in \cite{Mak2} and for convenience of the readers we
reproduce the proof of this theorem.

\proclaim{Theorem 10}
 Let $ R $ be a rational map.
Then: \roster \item    $\beta $ is an injection when restricted to $ HD(R)\times
T(J(R)), $ \item  if $ R $ is structurally stable, then $ \beta:HD(R)\times T(J(R))
\rightarrow H^1(R) $ is an isomorphism.
\endroster
\endproclaim
 \demo{Proof}
1). Let $ A(J(R)) \subset L_1(\overline{\Bbb C}) $ be the subspace of functions
which are holomorphic on $ F(R). $ Then $ A(J(R)) $ is a Banach space with the $ L_1
-$norm. Furthermore, let $ A(R) \subset A(J(R)) $ be the subspace of rational
functions. In other words, $ A(R) $ is the linear span of the functions $
\gamma_a(z) = \frac{a(a - 1)}{z(z - 1)(z - a)}, $ where $ a \in J(R). $ Then by the
Bers Density Theorem  $ A(R) $ is an everywhere dense subspace of  $ A(J(R)). $
\par Now let  $ \mu \neq 0 \in \ker(\beta)\cap \{HD(R)\times T(J(R))\}; $ then we have
$$
F_\mu(R(z)) = R'(z)F_\mu(z)
$$
and hence $ F_\mu = 0 $ on the set of all non-parabolic periodic points, and hence $
= 0 $ on the Julia set as well. Now if $ F(R) =\emptyset, $ then $ F_\mu = 0 $
identically on $ \overline{\Bbb C}; $ and using the fact that $ \mu =
\overline{\partial} F_\mu $ (in the sense of distributions) we have $ \mu = 0.$
\par If $ F(R) \neq \emptyset, $ then on $ A(J(R)) $ the functional
$ L_\mu(\phi) = \iint\mu\phi, $
 satisfies $ L_\mu(\gamma_a(z)) = F_\mu(a) = 0, $ for any $ a \in J(R). $
By the Bers density theorem we have $ L_\mu = 0 $ on $ A(J(R)), $ hence $ \mu = 0 $
almost everywhere on $ J(R). $ Hence we obtain $ \mu \in HD(R) $ represents the zero
- functional on the space $ A(\Omega(R)) \subset A(J(R)). $ By the discussion above,
we have $ \mu = 0.$
\par 2) If $ R $ is structurally stable, then $ \dim(HD(R)\times J_R) = \dim(H^1(R))
= 2deg(R) - 2. $ By 1) the
operator $ \beta $ is linear and injective, hence an isomorphism.
\enddemo

Now, let all critical points $ c_i $ be simple. Then there exists a decomposition $
\frac{1}{R'(z)} = \omega + \sum\frac{b_i}{z - c_i}, $ where $ \omega =
\frac{1}{R'(\infty)} $ is the multiplier of $ \infty $ and $ c_i $ are the critical
points (by the residue theorem $ b_i = \frac{1}{R''(c_i)}).$ For $ i = 1, ...,
2deg(R) - 2 $ let $ h_i(z)
 = \frac{1}{R'(z)} - \frac{b_i}{z - c_i}. $
\proclaim{Proposition 11} For any rational map $ R $ with simple
critical points, the following statements hold, \roster \item Let
$ \gamma_a(z) = \frac{a(a - 1)}{z(z - 1)(z - a)} \in
L_1({\overline{\Bbb C}}) $ where $ a \in {\Bbb C}\backslash \{0,
1\} $ is not a critical point. Then
$$ R^*(\gamma_a(z)) = \frac{\gamma_{R(a)}(z)}{R'(a)} +
\sum_ib_i\gamma_a(c_i)\gamma_{R(c_i)}(z).
 $$
\par Let $ \tau_a(z) = \frac{1}{z - a}, $   where $ a \in \Bbb C $ is not a critical point.
Then
$$
R^*(\tau_a(z)) = \frac{\tau_{R(a)}(z)}{R'(a)} + \sum_i
b_i\tau_a(c_i)\tau_{R(c_i)}(z).
$$
\item If $ a = c_i $ is a critical point, then
$$ R^*(\gamma_a(z)) = (h_i(a) + b_i\frac{2c_i - 1}{c_i(c_i - 1)})\gamma_{R(a)}(z) +
\sum_{j \neq i}b_j\gamma_a(c_j)\gamma_{R(c_j)}(z), $$ and
$$
R^*(\tau_a(z)) = h_i(a)\tau_{R(a)}(z) + \sum_{j \neq
i}b_j\tau_a(c_j)\tau_{R(c_j)}(z),
$$
where  $ h_i(a) + b_i\frac{2c_i - 1}{c_i(c_i - 1)} = \lim_{a\to c_i}
\left(\frac{1}{R'(a)} + b_i\gamma_a(c_i)\right). $
\endroster
\endproclaim
\demo{Proof} See lemma 5 in {\cite{Mak2}}
\enddemo

 From the proposition 11 we have
$$
\beta(\mu)(z) = F_\mu(R(z)) - R^\prime(z)\cdot F_\mu(z) = -R'(a)\sum_i
b_iF_\mu(R(c_i))\gamma_a(c_i).\tag{*}
$$

\proclaim{Remark 12} Proposition 11 gives a different set of
coordinates for the spaces $ H^1(R) $ and $ HD(R)\times J_R. $
Namely the formula $ * $ above describes the isomorphism $ \beta^*
: HD(R)\times J_R \rightarrow {\Bbb C}^{(2deg(R) - 2)} $ by
$$\beta^*(\mu) = (F_\mu(R(c_1)), ..., F_\mu(R(c_{2deg(R) - 2}))
$$
\endproclaim


\heading{\bf Beltrami Differentials on Julia set} \endheading
\par Here we discuss the space $ T(J(R)). $ Each element
$ \mu \in T(J(R))$ defines an $ invariant $ with respect to the
Ruelle operator functional $ L_{\mu} $ on the space $ A(R), $
which is continuous in the topology of $ A_2 $ (recall that $ A_i
= (A(R), \vert\cdot\vert_i)).$ Continuity of $ L_\mu $ in the
topology of the space $ A_1 $ is crucial as regards  the question
of non-triviality of $ \mu.$ Indeed we have the following lemma.
\proclaim{Lemma 13} Let $ \mu \in T(J(R)), $ then $ \mu = 0 $ if
and only if $ L_\mu $ is a continuous functional on $ A_1.$
\endproclaim
\demo{Proof} Let $ L_\mu $ be continuous on $ A_1 $. Then by the Bers' density
theorem $ L_\mu $ is continuous on $ A(\Omega). $  By  lemma 9, there exists an
element $ \psi \in A(S_R), $ such that the functional $ L_\mu(\alpha) =
\iint\limits_{\Bbb C}\alpha\lambda^{-2}\overline{\psi}, $ and hence $ F_\mu(a) =
L_\mu(\gamma_a) = \iint\limits_{\Bbb C }\gamma_a\lambda^{-2}\overline{\psi} =
F_{\lambda^{-2}\overline{\psi}}(a). $ Hence $ \beta(\mu)(a) =
\beta(\lambda^{-2}\overline{\psi})(a) $  for any $ a \in S. $ The set $ S $ is an
infinite subset of the plane, hence the rational functions $\beta(\mu) $ and $
\beta(\lambda^{-2}\overline{\psi}) $ are equal. This contradicts to the injectivity
of $\beta, $ and the lemma is proved.
\enddemo

\par Now we begin to consider the relationship between the continuity of $
L_\mu $ for $ \mu \in {\bold J}_R $ and certain properties of the
Ruelle operator $ R^\ast : A_2 \rightarrow A_2. $ Recall that the
operator $ R^\ast $ acts as a linear endomorphism of $ L_1({\bold
J}(R)) $ with unit norm. \proclaim{Proposition 14} Let $ R \in X $
be a rational map with simple critical points.
Assume that $ R $ is not the Latt{\'e}s map. Then \roster \item $
T(J(R)) = \emptyset $ if and only if the Ruelle operator $ R^\ast
: A_2 \rightarrow A_2 $ is mean ergodic, \item Assume in addition
that $ F(R) \neq \emptyset, $ and $ m(P(R)) = 0, $ then $ m(J(R))
= 0, $ if and only if the modulus of the Ruelle operator $ \vert
R_\ast\vert : L_1(J(R)) \rightarrow L_1(J(R)) $ is mean ergodic.
\endroster
\endproclaim
\demo{Proof} (1). If $ T(J(R)) = \emptyset, $ then the subspace $ (I - R^\ast)(A_2)
$ is everywhere dense in $ A_2, $ and by  item (2) of the Mean ergodicity lemma we
are done.
\par  Now suppose
that $ R^\ast $ is mean ergodic on $ A_2. $ Let $ \mu \neq 0 \in
T(J(R)), $ then there exists an element $ \gamma \in A_2 $ such
that $ \iint_{\overline{\Bbb C}} \mu\gamma \neq 0. $ Let $
\gamma_n = A(n,R)(\gamma) $ be  Cesaro averages, then by the
Theorem A the limit $ \lim_{n\to\infty}\gamma_n = 0 $ in the
strong topology on $A_2.$  Beside we have
$$ lim_{n\to\infty}\iint_{\overline{\Bbb C}}\mu\gamma_n = lim_{n\to\infty}\iint_{\overline{\Bbb
C}}\mu\gamma = \iint_{\overline{\Bbb C}}\mu\gamma \neq 0.
$$
The contradiction above complete the proof of the (1).

(2). If $ m(J(R)) = 0, $ then the space $ L_1(J(R)) = \{0\}, $ and we are done.

Now let $ \vert R_\ast\vert : L_1(J(R)) \rightarrow L_1(J(R)) $ be mean ergodic.

Assume that $ m(J(R) > 0. $ Let $ \phi > 0 \in L_1(J(R)) $ be any
element with $ \iint_{J(R)} \phi = 1. $ Then the mean ergodicity
of $ \vert R^*\vert $ implies that the Cezaro averages $ A(N,
\vert R^*\vert, \phi) $ converges to an element $ \psi $ so that
\roster\item $ \vert R^*\vert\psi = \psi $ and \item $ 1 =
\iint_{J(R)} \phi = \iint_{J(R)}A(N, \vert R^*\vert, \phi) =
\iint_{J(R)} \psi.$
\endroster
Now under assumption of the proposition and by the Lemma 7 above
we have that $ m(C) = 0.$ Hence $ J(R) = D $ modulo the Lebesgue
measure. If $ W \subset J(R) $ is wandering then $ \iint_W \psi =
0, $ and hence $ \psi = 0 $ on $ D = J(R) $ which is contradiction
with (2) above.
\enddemo


We will now  show that the topologies $ \Vert\cdot\Vert_1 $ and $
\Vert\cdot\Vert_2 $ are "mutually disjoint." Denote by $ X_i $ the
closure of the space $ \left(I - R^*\right)\left(A(R)\right) $ in
the spaces $ A_1 $ and $A_2 .$ \proclaim{Proposition 15} Let $ R $
be a rational map and $ dim(A(S_R))\geq 1. $ Then the following
conditions are equivalent. \roster \item the map $ i=id : A_1
\rightarrow  A_2 $ maps weakly convergent sequences onto weakly
convergent sequences. \item $ i(X_1) \supset X_2, $ \item the
Lebesgue measure of the Julia set is zero.
\endroster
\endproclaim
\demo{Proof} Condition (3) trivially implies  conditions (1) and (2).
\par Assume condition (1)  holds. Then the dual map $ i^*: A^*_2 \rightarrow  A^*_1
 $ is continuous in the $\ast-$weak topologies on $ A^*_1 $ and $ A^*_2. $ Hence for any
  $ \mu \in A^*_2 \subset L_\infty(J), $ there exists an element $ \nu \in A^*_1
   \subset L_\infty(F) $ such that $ \nu = i^*(\mu) $ and
$$ \iint_J\mu\gamma = \iint_F\nu\gamma.$$
Then for any $ \gamma \in A(R) $ we have $ \iint_{\Bbb C}\gamma (\mu - i^*(\mu)) =
0. $ Let $ F_\mu(z) $ and $ F_\nu(z) $ be potentials. Then $ F_{{\big\vert}J(R)} =
\left(F_\mu(z)- F_\nu(z)\right)_{{\big\vert}J(R)} = 0 $ and if $ m(J(R)) > 0 $ we
have $ F_{\overline{z}} = 0 $ almost everywhere on $ J(R), $ where $
F_{\overline{z}} $ is defined in the sense of distributions. Hence we deduce:
$$ \mu - i^*(\mu) = 0 $$
almost everywhere on $ J(R). $ Since $F(R)\cap J(R) =\emptyset, $ we have  $ \mu = 0
$ almost everywhere and  we conclude that $ A^*_2=\{0\}. $ Hence $A_2 =\{0\}, $
which gives  $ m\left(J(R)\right) = 0.$
\par Now assume (2). Then the hypothesis implies that any invariant continuous functional
on $ A_1 $ generates an invariant line field on the Julia set
contradicting the injectivity of the Bers map. By the assumption $
R $ always has non-trivial qc-deformation and hence, we conclude
that $ m\left(J(R)\right) =0. $
\enddemo
\proclaim{Proposition 16} Assume that $ \dim(A(S_R)) \geq 1 $ and
$ m\left(J(R)\right) > 0 $ for the given rational map $ R.$ Then
there exist no invariant line fields on the Julia set if and only
if $ i^{-1}(X_2) \supset X_1.$
\endproclaim
\demo{Proof} If there exist no invariant line fields, then $ X_2 = A_2.$ Now assume
$ i^{-1}(X_2) \supset X_1;$ then existence of an invariant line field would
contradict the injectivity of the Bers map.
\enddemo
\par We finish this chapter with the next theorem.

\proclaim{Theorem 17} Let $ R(z)  $ be a rational map and $c \in
J(R) $ be a critical point. Let $ S_L = \sum^{L}_{j =
0}\frac{1}{(R^{j})'(R(c))}. $ Assume that there exists a
subsequence $\{n_i\} $ of integers such that the sequence $
\{R^{n_i+ 1}(R(c))\} $ is bounded  and either \roster \item $
\lim_{i \rightarrow \infty} \vert (R^{n_i})'(R(c))\vert =\infty $
and
 $\overline{\lim}_{i \rightarrow \infty} \vert S_{n_i}\vert > 0 $ or
\item $ \vert (R^{n_i})'(R(c))\vert \sim C = Const $ for  $ i\to\infty $ and
$\overline{\lim}_{i \rightarrow \infty} \vert S_{n_i}\vert = \infty. $
\endroster
Then $ R $ is not  structurally stable (is an unstable map).
\endproclaim
\demo{Proof} Consider the one-dimensional family of deformations $ R_\lambda(z) =
R(z) + \lambda. $ Assume that $ R $ is  stable. Then there exist an $ \epsilon $ and
a holomorphic family $ h_\lambda : \overline{\Bbb C}\rightarrow \overline{\Bbb C} $
of qc-homeomorphisms such that for any $ \vert\lambda\vert < \epsilon, $
$$ R_\lambda(z) = h_\lambda\circ R\circ h_\lambda^{-1} (z).
$$
Let $ V(z) = \frac{\partial h_\lambda}{\partial\lambda}_{\vert\lambda = 0}(z); $
then the derivative of the equation above with respect to $ \lambda $ evaluated in $
\lambda = 0 $ gives the  equation
$$ V(R(z)) = 1 + R'(z)V(z).
$$
The function $ V(z) $ is continuous on $ \Bbb C $ and for any critical point $ c $
we have $ V(R(c) = 1. $
\par Now for any $ m $ and $ z \notin  R^{-m}(\infty), $ we calculate
$$ \frac{V(R^m(z))}{(R^m)'(z)} = V(z) + \frac{1}{R'(z)} + \frac{1}{(R^2)'(z)} + ... +
 \frac{1}{(R^m)'(z)}.
 $$
 Setting $ z = R(c) $ and $ m = n_i, $ we obtain in both cases a
 contradiction.
\enddemo

\heading{\bf Measures. Proof of theorem B.}
\endheading
\par Start again with a rational map $ R. $
Consider an element $ \gamma \in A(R) $ and the corresponding
Cesaro average sequence $ A_N(R)(\gamma)$ $=$ $\frac{1}{N}\sum_{i
= 0}^{N - 1} (R^*)^i(\gamma). $ Let $ C(U) $ be the space of
continuous functions defined on $ \overline U $ for a fixed
essential neighborhood $ U. $ Then any $*$-weak limit of $
A_N(R)(\gamma) $ on $ C(U) $ is called a {\it weak boundary} of $
\gamma $ respect to $ R^* $ over $ U; $ denote the set of all
limit measures  by $ \gamma(U,R).$ \proclaim{Proposition 18} Let $
R $ be a structurally stable rational map with non empty Fatou
set. Assume there exists a non-zero weak boundary $ \mu \in
\gamma(U, R^*) $ for an element $ \gamma \in A(R) $ and an
essential neighborhood $ U. $ Then the Lebesgue  measure $ m(J(R))
> 0 $ and there exists a non-trivial invariant line field on $ J(R).$
\endproclaim
\demo{Proof} Under the assumptions, there exists an essential $ U,
$  $ \gamma \in A(R) $ and  a subsequence $ N_i $ such that
\roster \item $ \iint\phi A_{N_i}(R)(\gamma) $ converges for any $
\phi \in C(U) $ and \item there exists $ \psi \in C(U) $ such that
$\lim_{i\rightarrow\infty} \iint\psi A_{N_i}(R)(\gamma) \neq 0. $
\endroster
    By  density of the space of compactly supported continuous function in the space $ C(U), $
we may assume that $ \psi $ has  compact support $ D \subset
\overline{U}. $ Extending  $ \psi $ to $ \overline{\Bbb C}
\backslash D $ by zero, we obtain $
\lim_{i\rightarrow\infty}\iint_{\overline{\Bbb C}}\psi
A_{N_i}(R)(\gamma) \neq 0. $ Hence the dual average $
A_N(B_{R})(\psi) = \frac{1}{N}\sum_{i = 0}^{N - 1} (B_{R})^i(\psi)
$ has non-zero $*$-weak limit element in the $ *$-weak topology on
$ L_\infty(J(R)). $ Let $ \mu \in L_\infty(J(R)) $ be this
non-zero limit element. Then $ \mu $ is fixed for $ B_{R} $ and $
\mu = 0 $ on $ F(R) $ by construction. Hence $ m(J(R)) > 0 $ and $
\mu $ defines the desired invariant line field.
\enddemo

It is not clear if the converse is true. We suggest the following
conjecture. \proclaim{Conjecture} Let $ R $ be a rational map with
non-empty Fatou set. The $ T(J(R)) = \emptyset $ if and only if
the weak boundaries $ \gamma(U,R^*) = 0 $ for all $ \gamma \in
A(R) $ and every essential neighborhood $ U. $
\endproclaim
\par In general the absence of invariant line fields on the Julia set implies mean
ergodicity of $ R^* $ on $ L_1(J(R)), $ and so it would be
interesting to understand the conditions implying  mean ergodicity
of $ R^* $ from the measure-theoretic point of view. To do this,
let us recall the definition of the following objects: \roster
\item $ U $ is an essential neighborhood of $ J(R) $ and \item $
H(U) $ consists of $ h\in C({\overline U}) $ such that $
\frac{\partial h}{\partial \overline{z}} $ (in sense of
distributions) belongs to $ L_\infty(U) $ \item $ H(U) $ inherits
the topology of $ C({\overline U}).$
\endroster
\par Measures $ \nu^i_l. $
\roster \item Let $ c_i $ and $ d_i $ be critical points and
critical values, respectively. Then  define $ \mu_n^i =
\frac{\partial}{\partial\overline{z}}((R^*)^n(\gamma_{d_i}(z)) $
(in sense of distributions).
\item Define by $ \nu^i_l $ the average $ \frac{1}{l}\sum_{k =
0}^{l - 1}\mu_k^i. $
\endroster
\subheading{Proof of Theorem B} Suppose that $ \nu^1_{l_k} $
converges in the $*$-weak topology on $ H(U) $ for
a subsequence $ \{l_k\} $ and an essential neighborhood $ U. $
Then the sequence of averages $ A_n(R)(\gamma_{d_1}(z)) \in L_1(U)
$ is weakly convergent. If $ m(J(R))
> 0, $ this means  $ A_n(R)(\gamma_{d_1}(z)) $ converges strongly in $ L_1(J(R)).
$ Let $ f = \lim_{n\to\infty}A_n(R)(\gamma_{d_1}(z), $ then by the
arguments of the theorem A either $ f = 0, $ or $ R $ is a
Latt{\'e}s. In the last case $ R $ is an instable map. Now let $
\mu \in HD(R)\times T(J(R)), $ then $ F_\mu(d_1) =
\iint_{\overline{\Bbb C}}\mu\gamma_{d_1}(z) =
\iint_{\overline{\Bbb C}}\mu A(n, R)(\gamma_{d_1}(z)) =
\lim_{n\to\infty}\iint_{\overline{\Bbb C}}\mu A(n,
R)(\gamma_{d_1}(z)) = 0. $ Contradiction with injectivity of the
Bers map.

\subheading{Proof of the Corollary B} Assume that the measures $
\nu^i_{l_k} $ converges in the $*$-weak topology on $ H(U) $ for
all $ i, $ a subsequence $ \{l_k\} $ and an essential neighborhood
$ U. $ Then the $ A_{l_k}(R)(\gamma_{d_i}(z)) $ converge strongly
in $ L_1(J(R). $

Now let $ \mu \neq 0 \in T(J(R). $ If  $ d \in J(R) $ is a
critical value, then by the arguments of the theorem C $ F_\mu(d)
= 0. $
Hence $ \beta(\mu) = 0, $ and $ \mu = 0. $ The contradiction with
assumption complete the proof.
\par Now assume $ T(J(R)) = 0. $ Let us show that
$ \nu_l^i \rightarrow 0 $ in the $ *$-weak topology on $ H(U) $
for any essential neighborhood $ U. $ Otherwise, there exists a
sequence $ \{l_k\}, $ an essential neighborhood $ U $ and a
function $ F \in H(U) $ such that
$$ \lim_{k\rightarrow\infty}\iint F\nu^{i_0}_{l_k} =
\lim_{k\rightarrow\infty}\iint
F_{\overline{z}}A_{l_k}(R)(\gamma_{d_{i_0}}(z)) \neq 0.
$$
Hence by the Mean Ergodicity Lemma  $
lim_{k\to\infty}A_{l_k}(R)(\gamma_{d_{i_0}}(z)) = \phi \neq 0. $
Then $ R $ is a Latt{\'e}s by the theorem A, and hence $ T(J(R)
\neq 0. $ The contradiction with assumption complete the proof of
the Corollary B.

\heading{\bf Proof Theorem C}
\endheading
 \par We begin by collecting some facts (see books of I.Kra
 "Automorphic forms and Kleinian Group" I. N. Vekua "Generalized
 analytic function.")
\proclaim{Facts} Denote by $ F_{\mu}(a) $ the following integral $
\iint_{\Bbb
 C}\mu(z)\tau_a(z)dzd{\overline{z}} $ where $\tau_a(z) = \frac{1}{z -
 a}$ for $ a \in {\Bbb C}$ and $ \mu \in L_{\infty}(J(R)). $
 Then
 \roster
 \item $ F_{\mu}(a) $ is a continuous function on $ {\Bbb C} $ and
 $\frac{\partial F_{\mu}(a)}{\partial {\overline{z}}} = \mu $ in the sense
 of distributions.
 \item $\mid F_{\mu}(a)\mid = O(\mid z\mid^{-1}) $ for large $ z.$
$\Vert F_{\mu}(a)\Vert_{\infty} \leq \Vert\mu\Vert_{\infty}M, $
where $ M $ does not depend on $ \mu $ and $ a \in {\Bbb C} $.
\item $\vert F_{\mu}(a_1) - F_{\mu}(a_2)\vert \leq
\Vert\mu\Vert_{\infty}C\vert a_1 - a_2\vert\vert\ln\vert a_1 -
a_2\vert\vert $, where $ C $ does not depend on $ \mu $ and $ a $.
\endroster
\endproclaim
\par  Denote by $ B  : L_\infty(J(R))
\to C(\Bbb C) $   the operator $ \mu \to F_{\mu}(a) $ and by $ X $
the image $ B(L_\infty(J(R)) $. Let $ W $ denote the space $ X $
with the following topology:
$$ \phi_n \to 0 \text{ iff } \Vert\phi_n\Vert_{\infty} \to 0
\text{ and } \frac{\partial\phi_n}{\partial\overline{z}} \to 0
\text{ in the $*$-weak topology of } L_\infty(J(R)).
$$
\proclaim{Lemma 19} \roster \item $ W $ is a complete locally
convex vector topological space. \item B is a compact operator
mapping $L_\infty(J(R))$ onto $ W.$
 \item Any bounded set $ U \subset W $ is
precompact.
\endroster
\endproclaim
\demo{Proof} Item (1) is obvious.
\par 2). Let $ U \subset L_{\infty}(J(R)) $ be bounded. Then $ U $ is
precompact in the $*$-weak topology on $ L_{\infty}(J(R)) $.
Furthermore, from item (2) of Facts, we have that $ B(U) $ forms a
uniformly bounded and equicontinuous family of continuous
functions. This means $ B(U) $ is precompact in the topology of
uniform convergence.
\par 3). Boundedness in $ W $ means in particular  that the set
$$ V = \{\frac{\partial\phi}{\partial\overline{z}} \text{ in sense of
distributions, for }\phi \in U \}
$$
forms a bounded set in the $*$-weak topology of $L_{\infty}(J(R)).
$ Hence $ V $ is bounded in the norm topology of $L_{\infty}(J(R))
$. We finish the  lemma by using item (2) and the fact that $ \phi
= B(\phi_{\overline{z}}).$
\enddemo
\par Define an operator $ T $ on $ X $ as follows
$$
T(F_\mu(a)) = F_{B_R(\mu)}(a) = \iint_{\Bbb C}B_R(\mu)\tau_a
=\iint_{\Bbb C}\mu R^*(\tau_a(z)).
$$
By Proposition 11 we have
$$
T(\phi) = \frac{\phi(R(a))}{R'(a)} - \sum\frac{b_i\phi(R(c_i))}{a
- c_i},
$$
where $ b_i $ are residue the function $ \frac{1}{R'(a)} $ in the
critical point $ c_i.$  For example for $ R(z) = z^2 +c $ we have
$ T(\phi)(a) = \frac{\phi(R(a)) - \phi(c)}{R'(a)}.$
\proclaim{Remark 20} By  definition it may be seen that
$$ \{T^n(\phi), n = 0, 1, ...\}
$$
forms bounded set in $ W. $
\endproclaim
\proclaim{Lemma 21} $ T $ is a continuous endomorphism of $ W $.
\endproclaim
\demo{Proof} Let $ F_{\mu_i} \to 0 $ in $ W; $ then $
\Vert\mu_i\Vert \leq C < \infty $ and hence $\{T(F_{\mu_i})\}$
forms  a precompact family in $ W. $  Let $ \psi_0 $ be a limit
point of this set. Then
$$ \psi_0(a) = \lim_j T(F_{\mu_{i_j}}) = \iint_{\Bbb
C}\mu_{i_j}R^*(\tau_a) \to 0 (*-\text{weak topology}).
$$
Thus $ \psi_0 = 0.$
\enddemo

Now let $ P $ be a strongly convergent polynomial. Let $ \Delta
\subset \overline{\Bbb C} $ be a closed topological disk centered
at infinity which does not contain critical point of any iteration
$ P^n. $  Then the family of functions
$$
 s_n(z) = \sum\frac{\vert b^n_i\vert}{\vert z - c^n_i\vert}
 $$
is uniformly bounded on $ \Delta, $ where  $ \sum\frac{b^n_i}{z -
c^n_i}  = \frac{1}{(P^n)'(z)}.$

 \proclaim{Lemma 22} Assume $ s_n(a)
\leq C < \infty $ for all $ n $ for a given polynomial $ P. $ Then
the Cesaro average $ A_N(\tau_a) $ converges with the $
L_1(J)-$norm.
\endproclaim
\demo{Proof} In the notations above, we have
$$\split
T(F_\mu)(y) &= \iint_J\frac{B(\mu)}{z - y} dz\land d{\overline{z}}
=
\iint_J\mu R^\ast(\tau_y)dz\land d{\overline{z}}\\
&= \frac{F_\mu(P)(y)}{P^\prime(y)} - \sum\frac{b_i
F_\mu(R(c_i))}{y - c_i} = \sum\frac{b_i(F_\mu(R(y)) -
F_\mu(R(c_i)))}{y - c_i}.\endsplit
$$
Now consider the sequence of functionals $ l_i(F) = (A_i(T)(F))(a)
$ on $ W. $ Under assumption of the lemma, we have
$$
\vert l_i(F)\vert \leq 2\frac{1}{i}\sum_{j=0}^{i -1}s_j(a)\sup_{w
\in {\Bbb C}} \vert F(w)\vert,
$$
 so the family of functionals $\{l_i\} $ can be extended onto the space $
C(\overline{\Bbb C}) $ of continuous functions on $ \overline{\Bbb
C} $ to the family of uniformly bounded functionals. Therefore we
can choose a subsequence $ l_{i_j} $ converging pointwise to some
continuous functional $ l_0. $ Note that $ l_0 $ is the fixed
point for the dual operator $ T^* $ acting on dual $ W^*.$ This
means that the sequence $ A_{i_j}(R^*)(\tau_a) $ weakly converges
in $ L_1(J), $ and hence by  the Mean ergodicity Lemma, the whole
sequence $ A_N(\tau_a) $ converges in norm to a fixed element of $
R^*. $
\enddemo


\par We are now ready to prove the  theorem C.
\subheading{Proof of Theorem C}

 It is enough to show convergence
of $ A_N(R^*)(\frac{a(a - 1)}{z(z - 1)(z - a)}) $ for any fixed $
a \in S.$

 Let us denote by $ Y $ the subset
of elements from $ L_1(J(P)) $ on which the averages $ A_N(P^\ast)
$ are convergent. Note that $ Y $ is a closed space such that
family $ A_N(P^\ast) $ forms equicontinuous family of operators.

\par We claim that {\it for any $ a \in S $ the elements
$\gamma_{a}(z) $ belong to $ Y. $} \demo{Proof of the
claim}Otherwise, there would exist a continuous functional $ L $
on $ L_1(J(R)) $ and $ a_0 \in S $ so that $ L(\gamma_{a_0}) \neq
0 $ and $ Y \subset ker(L).$ Note that $ L $ is an invariant
functional (i.e. $L(R^\ast(f)) = L(f)$) such that for any $ f \in
L_1(J(R)), $ the element $ f - R^\ast(f) $ belongs to $ Y. $ Let $
\nu \in L_\infty(J(R)) $ be the element corresponding to $ L; $
then $ \nu $ is a fixed  for Beltrami operator $ B_P $ and hence
the function $ F_\nu(a) =\iint\nu\tau_a $ is a fixed point for the
operator $ T $ i.e.
$$
\frac{F_\nu(P(a))}{P'(a)} - \sum\frac{b_iF_\nu(d_i)}{a - c_i} =
F_\nu(a).
$$
Let $ d \in \Delta $ be a point such that $ P(d) \in \Delta$; then
by the arguments above $ F_\nu(d) = F_\nu(P(d)) = 0 $. Therefore
the rational function $ \Phi(a) = \sum\frac{b_iF_\nu(d_i)}{a -
c_i} $ has a too large number of zeros. This immediately implies
$\Phi(a) \equiv 0 $ and the function $ F_\nu $ satisfies the
equation
$$ \frac{F_\nu(P(a))}{P'(a)} = F_\nu(a). $$
Finally we have that $ F_\nu $ is zero on the set of all repulsive periodic points,
hence on the Julia set,  hence everywhere because $ F_\nu $ is holomorphic on the
Fatou set. Thus $ 0 \neq L(\gamma_{a_0}) = (a_0 - 1)F_\nu(0) - a_0F_\nu(1) +
F_\nu(a_0) = 0, $ contradiction.
\enddemo

To complete the theorem C we need the following lemma:
\proclaim{Lemma 23} Let $ T(J(P)) \neq 0 $ for a strongly
convergent polynomial $ P. $ Then there exists a regular fixed
element for the Ruelle operator $ R^*.$
\endproclaim
\demo{Proof} Let $ \mu \neq 0 \in F(J(R)) $ be non-trivial
Beltrami differential, then there exists $ a \in {\Delta} $ such
that $ 0 \neq F_\mu(a) = \iint\mu\gamma_a(z). $ Then $ f =
lim_{n\to\infty} A_n(R^*)(\gamma_a(z)) \neq 0 $ is a fixed point
for the Ruelle operator. Let us show that  the total variation of
$ \overline{\partial}f $ is bounded. For this is sufficient to
show that the total variation of $ A_n(P^*)(\gamma_a(z)) $ is
bounded independently on $ n. $ Indeed we have,
$$
\vert \iint \phi \overline{\partial}P^*(\tau_a)\vert = \vert\iint
\overline{\partial}\phi P^*(\tau_a)\vert = \vert
\frac{\phi(P(z))}{P^\prime(a)} - \sum\frac{b_iF_\nu(d_i)}{a -
c_i}\vert \leq s_1(a)\Vert\phi\Vert,
$$
where $ \phi $ is any differentiable function. Hence  by the
induction we have desired result.

Now the contradiction with the theorem A complete proof of the
theorem C.
\enddemo
\par We will now give sufficient conditions on polynomial to be strongly convergent.
The conditions will be given in terms of the Poincar{\'e} series
of the rational map. We begin with the following calculations.
\proclaim{Lemma 24} Let $ R $ be a rational map with no critical
relations and simple critical points. Let $ c $ be a critical
point of $ R $ and $ d \in (R^k)^{-1}(c) $ be any point for some
fixed $ k. $ Then for any fixed $ m $ the coefficient $ b $
corresponding to the entry $ \frac{1}{z - d} $ in expression $
s_m(z) $ has the following type
$$
b = \frac{1}{(R^m)''(d)} = \frac{1}{(R''(c))(R^{m - k -
1})'(R(c))((R^k)'(d))^2}.
$$
\endproclaim
\demo{Proof} By a residue calculations at the point $ d. $
\enddemo
\par Let us recall the backward and forward Poincar{\'e} series for the given rational map $ R.$
\proclaim{Definition} Forward Poincar{\'e} series $ P(x, R) $
$$ P(x, R) = \sum_{n = 0}^\infty\frac{1}{\vert (R^n)'(R(x))\vert}. $$
Backward Poincar{\'e} series $ S(x, R). $

Let $ \vert R^*\vert = R_{1,1}$ be the modulus of the Ruelle
operator, then
$$ S(x, R) = \sum _{n = 1}^\infty \vert R^*\vert^n({\bold 1}_{\pmb{\Bbb C}})(x) = \sum_{n = 1}^\infty\sum_{R^n(y) = x}\frac{1}{\vert (R^n)'(y)\vert^2}.$$
\endproclaim
\subheading{Proof of Proposition C1} Let us again consider the
function $ s_n(a)=\sum\frac{\vert b_i\vert}{\vert a - c_i\vert} $
and let $ B_n = \sum\vert b_i\vert. $ Then by lemma above we have.
$$ B_n = \sum_{c \in Cr(R)}\frac{1}{\vert R''(c)\vert}\sum_{j = 1}^{n -1}\frac{1}{\vert (R^{n - j - 1})'(R(c))\vert}\sum_{R^j(y) = c}\frac{1}{\vert(R^j)'(y)\vert^2},
 $$
hence we have the following formal equality
$$ \sum_{n = 2} B_n = \sum_{c \in Cr(R)}\frac{1}{\vert R''(c)\vert}S(c, R)\otimes P(c,
R)
$$
where  $\otimes $ means Cauchy product of series.

\proclaim{Corollary C} Under condition of the theorem C above
assume that for any critical point $ c $ there exists a constant $
M_c $ so that
$$ \frac{1}{\vert (P^n)^\prime\vert} \leq \frac{M_c}{n}  \text{ and } \vert R^*\vert^n({\bold 1}_{\Bbb
C})(x) \leq \frac{M_c}{n},
$$
then $ P $ is a strongly convergent polynomial.
\endproclaim

\demo{Prof} The statement follows from properties of Cauchy
product; e.g. the Cauchy product of two harmonic series is
divergent but has uniformly bounded elements. We emphasize that
evidently there  is no rational maps for which the forward
Poincar{\'e} series is equivalent to harmonic series for any
critical point.
\enddemo
\subheading{Proof of Corollary C2}

\par Let $ R $ be a rational map and
$ c_i, i =1,..., 2deg(R) -2, $ and $ d_i, i =1,...,  2deg(R) -2, $ be critical
points and critical values, respectively, and let  $ z  =\infty $ be a fixed point
with multiplier $ \lambda. $ Then by induction we have.
$$ \frac{1}{R'(z)} = \lambda + \sum_i\frac{b_i}{z - c_i} =
\lambda + \sum_i\frac{1}{R''(c_i)}\frac{1}{z - c_i}, $$
$$ ... $$
$$ \frac{1}{(R^n)'(z)} = \lambda^n + \sum_i\sum_{k = 0}^{n - 1}
\left(\sum_{y \in R^{-k}(c_i)}\frac{1}{(R^{n})''(y)}\frac{1}{z - y}\right). $$
We finished the proof of the first equality with the following
lemma.
 \proclaim{Lemma 25} For any $ k < n $
$$
\multline
\sum_{y \in R^{-k}}\frac{1}{(R^{n})''(y)}\frac{1}{a - y} =
\frac{1}{R''(c_i)}\frac{1}{(R^{n - k - 1})'(d_i)}\sum_j\frac{(J_j')^2(c_i)}{a - J_j(c_i)} =\\ =\frac{1}{R''(c_i)}\frac{1}{(R^{n - k - 1})'(d_i)}(R^*)^k(-\tau_a)(c_i),
\endmultline
$$
where the $ J_j $ are branches of $ R^{-k}, \tau_a(z) = \frac{1}{z - a} $ and $ R^*
$ is Ruelle operator.
\endproclaim
\demo{Proof} Lemma 24 and the above equalities.
\enddemo

\par We now prove the second equality. By proposition 11, we  calculate as follows:
$$
\multline
\left(R^*\right)^0(\tau_a)(z) =\tau_a(z), \left(R^*\right)(\tau_a)(z) = \frac{1}{R'(a)(z - R(a)} - \sum_i\frac{b_i}{(a - c_i)(z - R(c_i))}\\
\left(R^*\right)^2(\tau_a)(z) = \frac{1}{(R^2)'(a)(z - R^2(a)} - \frac{1}{R'(a)}\sum_i\frac{b_i}{(R(a) - c_i)(z - R(c_i)} - \\
\sum_i\frac{b_i}{(a - c_i)}R^*(\frac{1}{(z - R(c_i))}
\endmultline
$$
and by induction,
$$\multline
\left((R^*)\right)^n(\tau_a)(z) = \frac{1}{(R^n)'(a)(z - R^n(a)} -\\
-\sum_ib_i\left(\frac{1}{(R^{n -1})'(a)(R^{n - 1}(a) - c_i)(z - R(c_i))} +
  ... + \frac{1}{a - c_i}(R^{n - 1})^*\left(\frac{1}{z - R(c_i)}\right)\right).
\endmultline
$$
Then summation with respect to $n $ gives the desired equality.


\heading{\bf Few Open Questions on spectrum of Ruelle and Beltrami
operators and so on}\endheading

The next proposition is completely trivial and we are left proof
to the reader.

\proclaim{Proposition 26} The spectrum of the operator $ R^* :
L_1(\overline{\Bbb C})\to L_1(\overline{\Bbb C}) $ is the closed
unit disk $ \Delta. $ Any interior point $\lambda \in \Delta $ is
the eigenvalue.
\endproclaim

Now let $ Y = \text{linear span} \{\gamma_a(z), a \in P(R)\},
\subset L_1(\overline{\Bbb C}) $ and $ Z $ be closure of $ Y, $
then we have the following connection of the dynamic of $ R $ and
spectrum of its Ruelle operator.

\proclaim{Proposition 27} \roster \item $ R^* $ maps $ Z $ into $
Z. $ \item The spectrum of $ R^*: Z\to Z $ is finite pure point
spectrum if and only if  $ R $ is postcritically finite map.
\endroster
\endproclaim
\demo{Proof} (1). This case follows immediately from the
proposition 11.

(2) If $ R $ is a postcritically finite, then $dim(Z) < \infty $
and we are done.

Assume That $ P(R) $ is infinite set, then there exists a point $
a \in P(R) $ with infinite forward orbit and hence $ \gamma_n =
R^*(\gamma_a(z)) $ are different pairwise. Let $ \lambda_i, i = 1,
..., N $ be the constants so that $ \gamma_a(z) = \sum_i
\lambda_i\phi_i, $ where  $ \phi_i $ are eigenfunctions. Now let $
\alpha_i, i = 1, ..., N $ be the points of the spectrum, here $
\alpha_0 = 0 $ if zero belongs to the spectrum. Then we have
infinitely many equations:
$$ \gamma_n = \sum_i\alpha^n\lambda_i\phi_i
$$
Hence there exists an $ i_0 \neq 0 $ so that the eigenfunction $
\phi_{i_0} $ is rational and indeed is a linear combination of
finite number of $ \gamma_k.$  Moreover $R^*(\phi_i) =
\alpha_{i_0}\phi_{i_0}. $ If $ A $ be the set of the poles of $
\phi_{i_0}, $ then $ A\cap \{\cup_n R^n(a)\} \neq\emptyset. $ By
the proposition 7 we have $ R^n(A) \subset A\cup \{\text{critical
values of} R^n\} $ for any $ n
> 0. $ Hence the forward orbit of $ a $ is finite which is a
contradiction.
\enddemo

It is interesting describe the spectrum of the operator $ R^*:
Z\to Z. $

Is it true that spectrum (not pure pointed) is finite only in a
case of postcritically finite map ?

Is it possible that for a rational map $ R $ the operator $ R^*:
Z\to Z $ is a compact operator with infinite spectrum?

\subheading{Spectrum of Beltrami operator}

Let us recall that a rational map $ R \in Rat_d $ is in {\it
general position} iff the cardinality of its critical values $ =
2d-2. $ Particularly that means that any two critical points of $
R $ have different images.

\proclaim{Definition} Let $ R \in Rat_d, $ then the Hurwits class
$ H(R) $ of the map $ R $ is the following space
$$
H(R) = \{g \in Rat_d, {\text{ there exist homeomorphisms }} \phi,
\psi, R\circ\phi = \psi\circ g\}.
$$
\endproclaim

Note that the Hurwitz class was introduced in holomorphic dynamic
by A.Eremenko and M.Lyubich \cite{EL} for entire function. We call
this set as a Hurwitz class because the Hurwitz theorem describes
classes of brunched coverings like this (see for example \cite B).

\proclaim{Lemma 28}
\roster\item Let $ R $ and $ Q $ be two
rational maps in general position, then $ Q \in H(R) $ if and only
if $ deg(R) = deg(Q). $ \item The set of the maps in general
position the  degree $ d $ forms open and everywhere dense subset
in $ Rat_d. $
\endroster
\endproclaim
\demo{Proof} See {\cite B}
\enddemo

\proclaim{Proposition 29} The following cases are equivalent
\roster \item The rational map $ R $ is hyperbolic.\item There
exists $ N > 1 $ such that $ R^N $ is $J-$stable in the space $
Rat_{d^N}.$ \item For any $ N > 1 $ the map $ R^N $ is $J-$stable
in the space $ Rat_{d^N}.$
\endroster
\endproclaim
\demo{Proof} The proof of this proposition is extremely trivial.
The (1) implies (2) and (3).

Now let $ N > 1 $ be the number from the conditions (1) or (2) of
the proposition. By the lemma above there exists a map $ Q \in
Rat_{d^N} $ which is in the general position and a homeomorphism $
\phi: J(R) \to {\Bbb C}, $ so that $ Q = \phi\circ R\circ
\phi^{-1}. $  Assume that there exists a critical point $ c \in
J(R), $ then only one possibility can occur $ R^{-N}(c) = c $ and
$deg(R) = 2.$ That means $ c $ is periodic superattractive point
and cannot belong to the Julia set. Contradiction with assumption.
\enddemo

Note that in non-hyperbolic case $J-$ stability in $ Rat_d $ means
that the dimension of the space of fixed points for Beltrami
operator should be equal to the number of critical points on the
Julia set and this last number growthes exponentially with respect
to iterations. Indeed the dimension of the space of fixed points
for Beltrami operator is equal to the number of {\it critical
values} on and this number growthes linearly with respect to
iterations.

Probably the following assumption: $ R^n $ is $ J-$stable in $
H(R^n) $ for any $ n > 1 $ implies the hyperbolicity  of $ R. $

Now let $ R $ be a rational map, assume that $ R^k $ is $J-$
stable in $ H(R^k),$ is it true that $ R^n $ is $ J-$stable in $
H(R^n),$ for $ n > k $? What is possible to say in the case $ n <
k $?

\proclaim{Remark 30} Note that the question above are also the
questions about point spectrum of Beltrami operator, that is if $
R^k $ is $J-$ stable in $ H(R^k),$ then the Beltrami operator $
B_R $ has an eigenvalue which is $k-$root of unity.
\endproclaim

\enddocument